\documentclass[preprint]{imsart}

\RequirePackage{amsthm,amsmath}
\RequirePackage[round]{natbib}
\RequirePackage[colorlinks,citecolor=blue,urlcolor=blue]{hyperref}
\usepackage{amssymb}
\usepackage{graphicx}
\usepackage{amsfonts}
\usepackage{amsbsy}
\usepackage{mathrsfs}
\usepackage{lscape}
\usepackage{float}
\theoremstyle{plain}

\usepackage{multirow}
\usepackage{tablefootnote}
\usepackage{threeparttable}
\usepackage{pdflscape}
\usepackage{enumerate}
\usepackage{pdfpages}
\usepackage{soul}
\usepackage[shortlabels]{enumitem}
\usepackage[french,english]{babel}
\usepackage{xcolor,colortbl}

\startlocaldefs
\numberwithin{equation}{section}
\theoremstyle{plain}
\newtheorem{theorem}{Theorem}
\newtheorem{definition}{Definition}
\usepackage[english]{babel}

\endlocaldefs

\newcommand{\bx}{\mathbf{x}}
\newcommand{\bX}{\mathbf{X}}

\newcommand{\bz}{\mathbf{z}}
\newcommand{\bs}{\mathbf{s}}
\newcommand{\bk}{\mathbf{k}}

\newcommand{\bH}{\mathbf{H}}
\newcommand{\bV}{\mathbf{V}}
\newcommand{\bbeta}{\boldsymbol{\beta}}
\newcommand{\bpsi}{\boldsymbol{\psi}}

\newcommand{\bSigma}{\boldsymbol{\Sigma}}
\renewcommand{\P}{\mathrm{P}}
\newcommand{\EE}{\mathbb{E}}
\newcommand{\RR}{\mathbb{R}}
\newcommand{\Var}{\mathrm{Var}}
\newcommand{\dd}{\mathrm{d}}
\newcommand{\yes}{yes }
\renewcommand{\no}{no}

\begin{document}


\begin{frontmatter}

\title{Regularization techniques for inhomogeneous (spatial) point processes intensity and conditional intensity estimation}
\runtitle{Regularization techniques for spatial point processes}

\begin{aug}
\author[A]{\fnms{Jean-François} \snm{Coeurjolly} \thanks{JF Coeurjolly is supported by ANR Labex Persyval-lab.}\thanks{JF Coeurjolly thanks all the organizers of the "Journées MAS 2022" and the members of the group MAS from SMAI for the opportunity offered to us to prepare this short review/tutorial paper.}},
\author[B]{\fnms{Isma\"ila} \snm{Ba}}
\and
\author[C]{\fnms{Achmad} \snm{Choiruddin}}
 
 \address[A]{Université Grenoble Alpes, LJK, 38000 Grenoble, France}

\address[B]{Department of Mathematics and Statistics, York University, Toronto, Ontario, Canada}

\address[C]{Department of Statistics, Institut Teknologi Sepuluh Nopember, 60111 Surabaya, Indonesia}




\end{aug}

\selectlanguage{english} 
\begin{abstract}  
\; Point processes are stochastic models generating interacting points or events in time, space, etc. Among characteristics of these models, first-order intensity and conditional intensity functions are often considered. We focus on inhomogeneous parametric forms of these functions assumed to depend on a certain number of spatial covariates. When this number of covariates is large, we are faced with a high-dimensional problem. This paper provides an overview of these questions and existing solutions based on regularizations.
\end{abstract}
\selectlanguage{french} 
\begin{abstract}
Les processus ponctuels constituent une classe de modèles stochastiques permettant de modéliser des évènements dans le temps, l'espace, etc en interaction. Parmi les caractéristiques d'un processus ponctuel, l'intensité et l'intensité conditionnelle d'ordre un sont souvent considérées. Nous nous concentrons ici sur des formes paramétriques inhomogènes de ces fonctions que nous supposons dépendre d'un certain nombre de covariables spatiales. Lorsque ce nombre est élevé, nous faisons face à un problème de grande dimension. Ce papier a pour objectif de présenter un aperçu de ces problèmes et solutions existantes.
\end{abstract}

\end{frontmatter}
\selectlanguage{english}


\section{Introduction}
Spatial point processes are stochastic processes which model point patterns distributed
in a space say $S$ (usually a subset of $\mathbb R^d$), such as locations of crime events,  species of trees, earthquake occurrences, disease
cases, etc (see e.g.\cite{baddeley2015spatial,illian2008statistical}). Modeling and inferring the intensity or conditional intensity of a spatial point process often constitutes the first and important task in the description and analysis of a spatial point pattern \cite{coeurjolly2019understanding}.
Roughly speaking, the intensity function measures the probability of an event to occur at a specific location, say $u \in S$, while the conditional intensity measures the probability to observe a point at $u$ given $\bx$ the observed set of events (i.e. points). This paper is focused on inhomogeneous models and in particular on parametric (conditional) intensity models for which intensities can be explained by spatial covariates (e.g. altitude map, soil nutrients, etc if one is interested in modeling locations of trees in a forest, see Figure~\ref{fig:bei}). In recent years, observing a large number of spatial covariates has become easier and provides more information on the point pattern which is analysed. A particular issue in fitting a parametric model to the intensity/conditional intensity arises when the number of covariates is large. To overcome such an issue, methodology based on regularization has been developed for spatial point process intensity/conditional intensity modeling. This research covers development on the methodology, theory, and computation \cite{thurman2014variable,thurman2015regularized,choiruddin2018convex,daniel2018penalized,choiruddin2023adaptive,ba2023inference,choiruddin2023combining}. 

In the literature (e.g. \cite{moller2003statistical,baddeley2015spatial}), estimating the intensity or conditional intensity function are two questions which are, the more often, treated separately. This was mainly justified by the fact that the object to model, mathematical tools (Campbell or Georgii-Nguyen-Zessin equations, see~\eqref{eq:Campbell}-\eqref{eq:GNZ}), statistical methodologies and proofs appear, at first glance, really different. The objective of this paper is to make a short overview of regularization techniques applied to point processes by trying, as far as possible, to present the problems of estimating the two kind of intensity functions in a similar way in order to shed the light on their similarities. To present methodologies, setting, etc, some choices are made. For instance, even if more general regularization techniques and asymptotic regimes were considered in \cite{choiruddin2018convex,choiruddin2023adaptive,ba2023inference}, we present only a part of asymptotic results and for one type of penalty, namely for the adaptive lasso. However, we propose an extension of a consistency result for which we  consider a general asymptotic regime and stochastic regularization parameters.
The rest of the paper is organized as follows. Background on spatial point process is described in Section~\ref{sec:spp}. We detail the statistical inference, theoretical results, and numerical aspects in Sections~\ref{sec:low}-\ref{sec:high}.


\section{Spatial point processes} \label{sec:spp}

\subsection{Notation and intensity functions}

We consider spatial point processes in $\RR^d$. For ease of exposition, we view a point process as a random locally finite subset $\bX$ of a Borel set $S \subseteq\RR^d$, $d\geq 1$. For readers interested in measure theoretical details, we refer to e.g.\ 
\cite{moller2003statistical,daley2007introduction} or \cite{dereudre2019introduction}. This setting implies the following facts. First, we consider simple point processes (two points cannot occur at the same location). Second, we exclude manifold-valued point processes (like circular or spherical point processes), and marked point processes, even if most of the concepts and methodologies presented hereafter exist or can be straightforwardly adapted to such contexts.

Thus, $\mathbf X\cap B$ stands for the restriction of $\mathbf X$ to a set $B\subseteq S$ and we let $|B|$ denote the volume of any bounded $B\subset S$. Local finiteness of $\mathbf X$ means that 
$\mathbf X \cap B$ is finite almost surely (a.s.), that is the number of points $N(B)$ of $\bX \cap B$ is finite a.s., whenever
$B$ is bounded. We let ${\cal N}$ stand for the state space consisting of the locally finite subsets (or point configurations) of $S$. 

The distribution of $\bX$ can be characterized by the finite-dimensional distributions of counting variables, or by the void probability, i.e. the probability to have no point in any compact set. However, these are usually not accessible and it is easier to summarize (and estimate) interpretable statistical measures such as intensity functions and conditional intensity functions. A more rigorous introduction on intensities, Palm intensities and conditional intensities and their links with reduced moment measures, Palm measures and reduced Campbell measures can be found in~\cite{coeurjolly2019understanding}. To get quicker to the core of the paper, we introduce them through Campbell theorem and Georgii-Nguyen-Zessin formula which may be viewed as integrals characterizations.

\begin{theorem}[Campbell theorem and GNZ formula]
The $k$-th order intensity function $\rho^{(k)}$ and the $k$-th order Papangelou conditional intensity function $\lambda^{(k)}$ are defined such that for any measurable function $h^{(k)}: (\RR^d)^k \to \RR^+$ and $\tilde h^{(k)}: \mathcal N \times (\RR^d)^k \to \RR^+$, we have respectively
\begin{align}
\EE \Bigg\{ 
\sum_{u_1,\dots,u_k}^{\neq} &h^{(k)}(u_1,\dots,u_k)
\Bigg\} 
= \nonumber\\
&\int_{\RR^d} \dots \int_{\RR^d} h^{(k)}(u_1,\dots,u_k) 
\rho^{(k)}(u_1,\dots,u_k) \dd u_1 \dots \dd u_k. \label{eq:Campbell}\\
\EE \Bigg\{ 
\sum_{u_1,\dots,u_k}^{\neq} &\tilde{h}^{(k)}(\{u_1,\dots,u_k\},\bX\setminus \{u_1,\dots,u_k\})
\Bigg\} 
= \nonumber \\
&\EE \left\{
\int_{\RR^d} \dots \int_{\RR^d} \tilde{h}^{(k)}(u_1,\dots,u_k) \lambda^{(k)}(\{u_1,\dots,u_k\},\bX) \dd u_1 \dots \dd u_k \right\}. \label{eq:GNZ}
\end{align}
\end{theorem}

\begin{figure}[htbp]
\centering
\includegraphics[width=.45\textwidth]{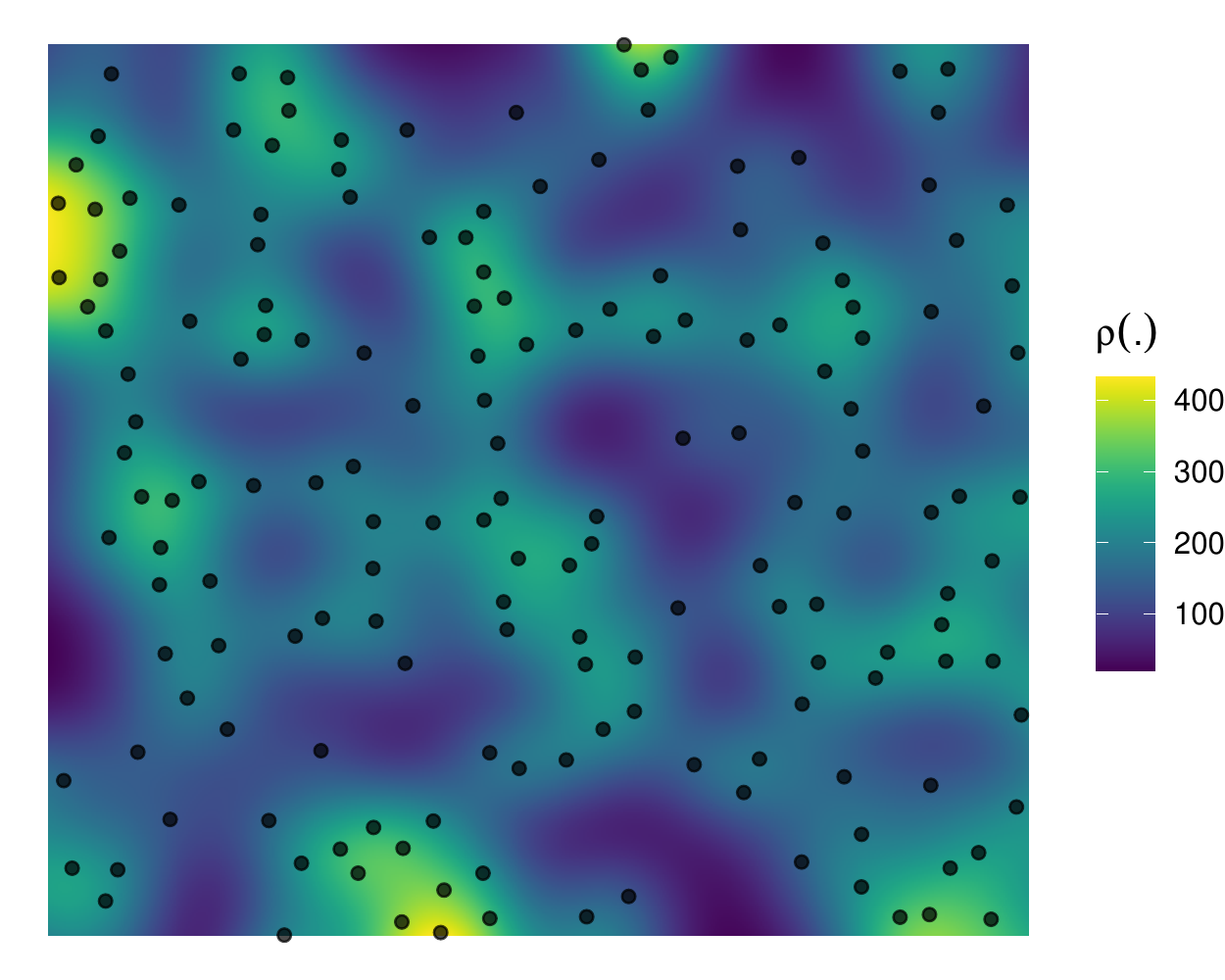}
\includegraphics[width=.45\textwidth]{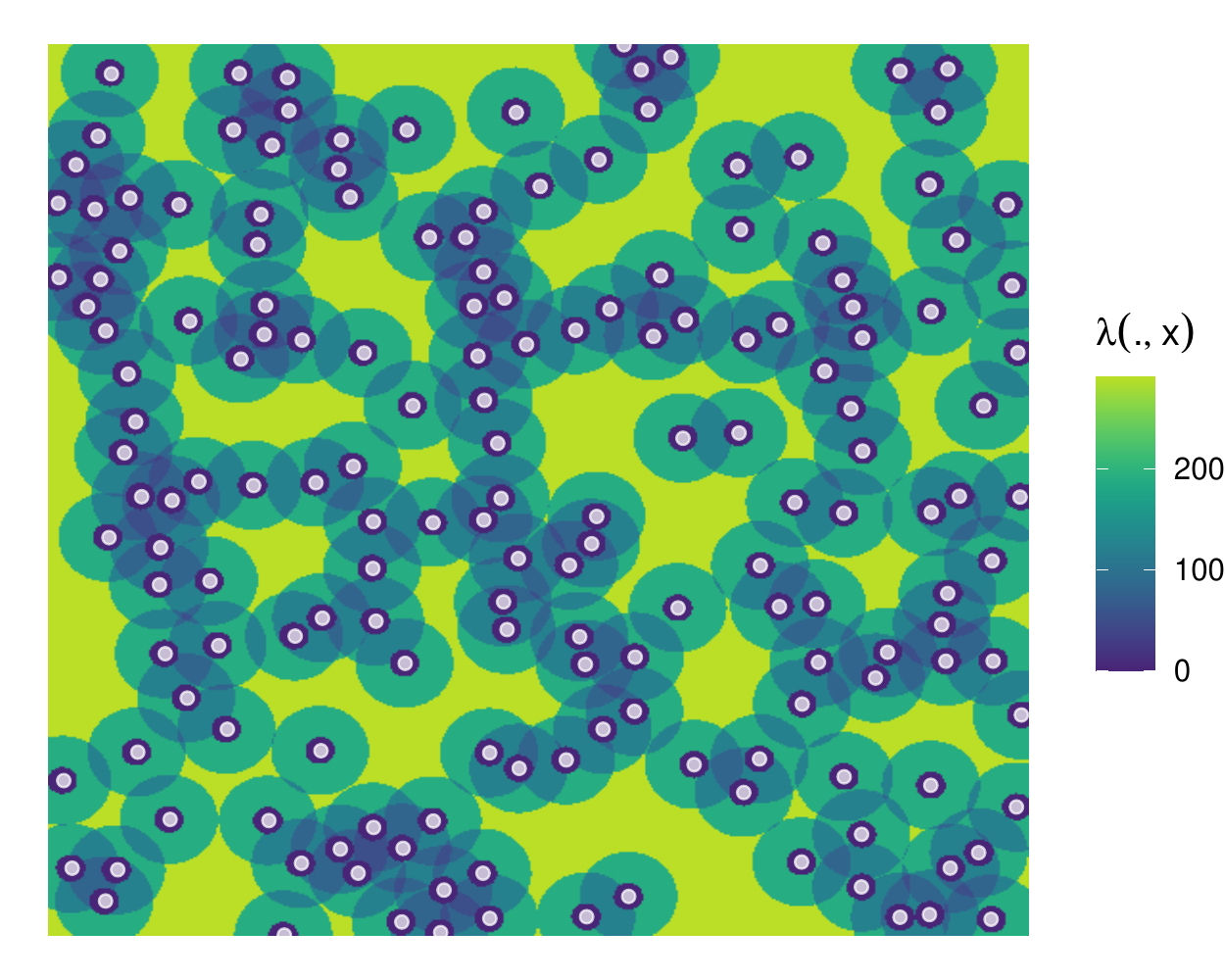}
\label{fig:int}
\caption{Left: Illustration of the concept of intensity function. The latent image is the true intensity function of the model (not specified) generating the pattern. The higher the intensity function the more likely a point can occur at this location; Right: Illustration of the conditional intensity function. The latent image depicts Papangelou conditional intensity function for a Gibbs model (not specified) given the configuration of points. The higher the intensity the more likely a point can be added at this location. We can observe that this model produces repulsive patterns as the conditional intensity is zero around points and quite small in balls around observed points.}
\end{figure}

When $k=1$, we more simply speak of the intensity function or the Papangelou conditional intensity function. It is relevant to have the following interpretation of such functions: $\rho(u)$ (resp. $\lambda(u,\bX)$) can be interpreted as the probability to observe a point in $B(u,\dd u)$ an infinitesimal ball centered at $u$ (resp. one point in $B(u,\dd u)$ given the rest of the configuration of points outside the ball is $\bX$). Similar interpretations are available when $k>1$. Equations~\eqref{eq:Campbell}-\eqref{eq:GNZ} can be combined to show that $\rho(u)=\EE\{\lambda(u,\bX)\}$ (also valid when $k>1$). Figure~\ref{fig:int} illustrates briefly the functions $\rho$ and $\lambda$. 

\subsection{Classes of models}

The reference model is the Poisson point process often defined as follows.
\begin{definition}\label{def:poisson}
Let $\rho$ be a locally integrable function on $\RR^d$. A point process $\bX$ satisfying the following statements is called the Poisson point process on $\RR^d$ with intensity function $\rho$:
\begin{itemize}
	\item for any $m\geq 1$, and for any disjoint and bounded $B_1,\dots,B_m  \subset S$, the random variables $\bX \cap {B_1},\dots,\bX \cap {B_m}$ are independent;
	\item $N(B)$ follows a Poisson distribution with parameter $\int_B \rho(u) \dd u$ for any bounded $B \subset S$.
\end{itemize}
\end{definition}
Poisson point processes model (eventually inhomogeneous) patterns with no interaction between points. As a consequence, it can be  easily proved that for such processes, $\rho^{(k)}(u_1,\dots,u_k)=\lambda^{(k)}(\{u_1,\dots,u_k\},\bX)=\prod_i \rho(u_i)$. Large classes of models exist to introduce dependence between points. A survey can be found in~\cite{coeurjolly2019understanding} and the references therein. This is debatable but, to our point of view, the main classes are: Cox processes (which includes Neymann-Scott, shot noise Cox or log-Gaussian Cox processes, see~e.g. \cite{moller2003statistical}) defined as Poisson point processes with random driven intensity obtained from a random field; Gibbs point processes (see~\cite{dereudre2019introduction}), which are (in a bounded domain) defined via a density with respect to a Poisson point process with intensity 1; Determinantal point processes (e.g.~\cite{lavancier2015determinantal}), for which intensities are defined through the determinant of a kernel function. We do not intend to define rigorously these models (see~e.g.~\cite{coeurjolly2019understanding} and references therein), however it is worth pointing out that these models are very different by the kind of interaction they model, their flexibility and, in particular as regards the concern of this paper, the fact that the intensity function and/or the Papangelou conditional intensity function is explicit or not. Table~\ref{tab:models} is an attempt to present the diversity of these models, their richness which also makes this research area attractive and fruitful.

\begin{table}[htbp]
\centering
\begin{tabular}{llll}
\hline
\rowcolor{blue!10} Model &  Type of interaction&Is $\rho(\cdot)$ explicit? &Is $\lambda(u,\bx)$ explicit?\\
\rowcolor{blue!10}
&& & \\
\hline
Poisson& no interaction & \yes & \yes \\
&&&\\
\hline
\rowcolor{black!10}Cox & attraction & \yes & \no \\
\rowcolor{black!10}&&&\\
\hline
Gibbs & attraction/repulsion & \no & \yes \\
&&&\\
\hline
\rowcolor{black!10}DPP &  repulsion& \yes & \yes and \no\\ 
\rowcolor{black!10}&&&\\
\hline
\end{tabular} 
\label{tab:models}
\caption{Attempt to classify most of spatial point process models. Column type of interaction refers to the type of patterns the corresponding model can produce. Last two columns answer the question of tractability of the intensity and conditional intensity function. "No" means that for most of models there is no explicit expression. The "yes and no" is more in between "yes" and "no". We can obtain an explicit expression but it quite complex to exploit.}
\end{table}

\subsection{Inhomogeneous parametric models}

\begin{figure}[htbp]
\centering
\includegraphics[width=.5\textwidth]{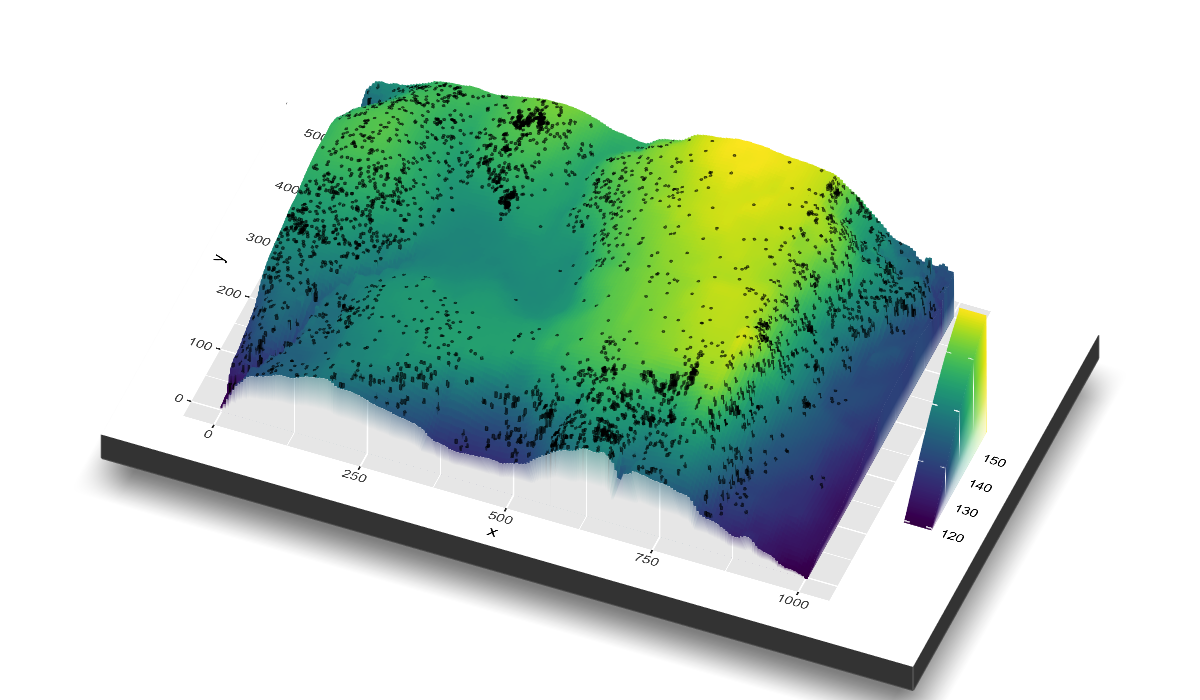}
\includegraphics[width=.49\textwidth]{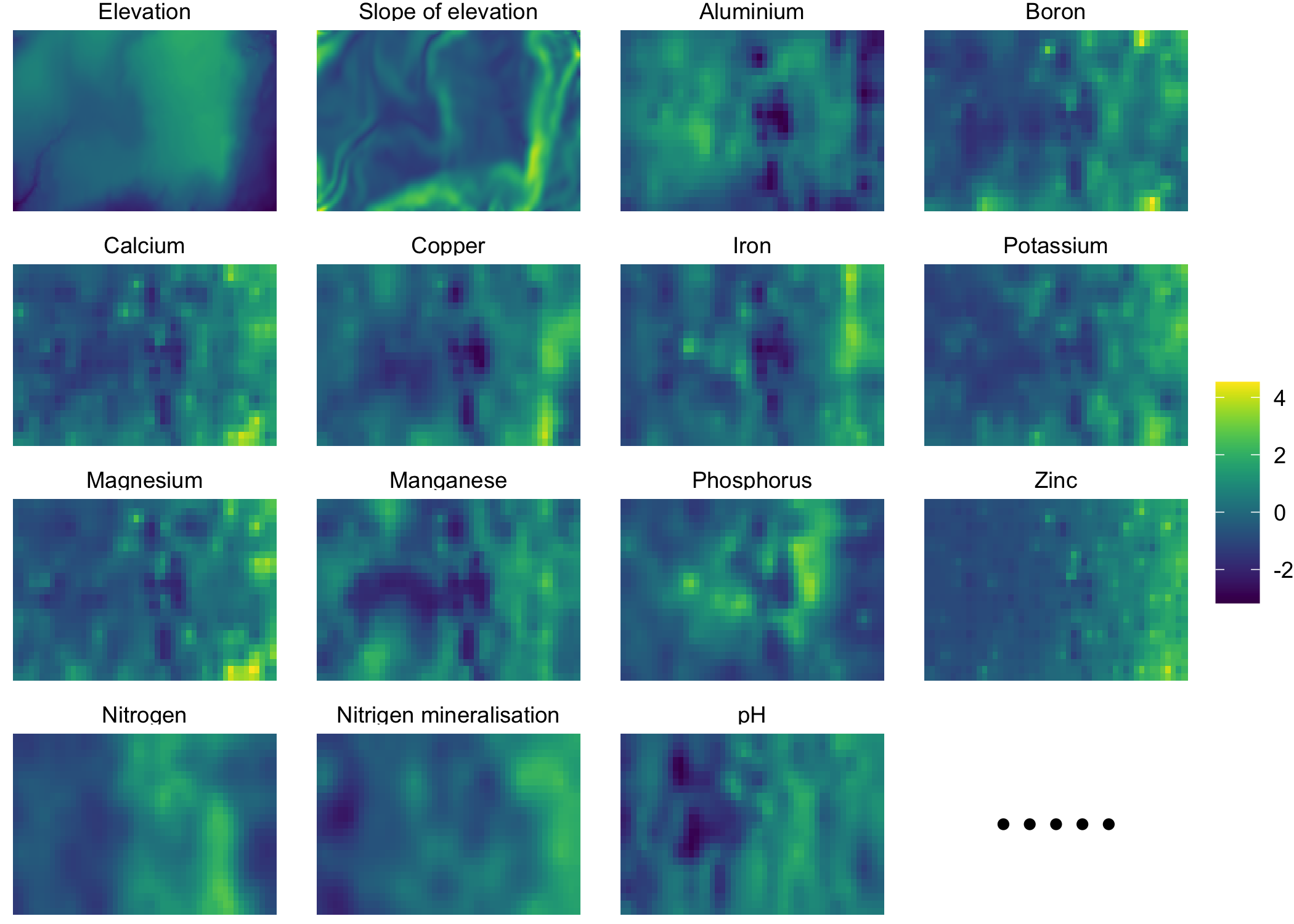}
\caption{Left: locations of one species of trees in Barro Colorado Island (dataset obtained from the Center of Tropical Forest) represented with the ground topography; Right: some of available spatial covariates (observed on the same observation domain) such as the altitude, slope of elevation and values of soil nutrients (Aluminium, Boron, etc).}
\label{fig:bei}
\end{figure}

Let us consider Figure~\ref{fig:bei} to motivate this section and paper. It is often the case in spatial statistics, that we observe a point pattern, here the locations of 3605 trees in a tropical forest (see~\cite{baddeley2015spatial} for more details on this dataset) together with spatial covariates which are information available on the whole observation domain. A quick look at Figure~\ref{fig:bei} is enough to be convinced of the inhomogeneity (and maybe non independence) characteristic of the point pattern and that it makes completely sense to relate the distribution of trees with covariates such as the elevation, the slope of elevation or levels of soil nutrients. 

In this application, we could be interested to model either the intensity and/or the Papangelou conditional intensity. We focus in this paper on exponential family models, where for any $u\in \RR^d$ and $\bx \in \mathcal N$ 
\begin{align}
\rho(u) &= \exp \left\{ \bbeta^\top \bz(u) \right\} 
\qquad \text{ and } \qquad
\lambda(u,\bx)  =\exp\left\{ \bbeta^\top \bz(u) + \bpsi^\top \bs(u,\bx) \right\}    \label{eq:model}
\end{align}
In both definitions, $\bbeta=\{\beta_1(u),\dots, \beta_p(u)\}^\top \in \mathbb R^p$ represents the main vector parameter of interest, \linebreak$\bz(u)=\{z_1(u),\dots, z_p(u)\}^\top, \; z_i:\RR^d \to \RR$ corresponds to the spatial covariates. We let $\bpsi \in \RR^l$ and $\bs(u,\bx) =\{ s_1(u,\bx),\cdots,s_l(u,\bx)\}^\top$ denote the parameter and the sufficient statistics defining the interaction term (more details are given below). To rewrite~\eqref{eq:model} into the same formalism we suggest the reformulation
\begin{align}
\rho(u;\bbeta_\rho) &= \exp \left\{ \bbeta_{\rho}^\top \bz_\rho(u) \right\} 
\qquad \text{ and } \qquad
\lambda(u,\bx;\bbeta_\lambda)  =\exp\left\{ \bbeta_\lambda^\top \bz_\lambda(u,\bx) \right\}    \label{eq:model2}
\end{align}
where $\bbeta_\rho \in \RR^p$ (resp. $\bbeta_\lambda \in \RR^{p+l}$), $\bz_\rho(u)=\bz(u)$ and  $\bz_\lambda(u,\bx)=\{\bz(u)^\top, \bs(u,\bx)^\top\}^\top$. In the rest of the paper, each time we index a vector, matrix, random quantity by $\rho$ (resp. $\lambda$) means that we refer to the estimation of $\rho$ (resp. $\lambda$). And when a comment applies to the two problems, we write $\bullet$. Hence, for instance $\bbeta_\bullet$ stands either for $\bbeta_\rho$ or $\bbeta_\lambda$. We point out that non-exponential family models can be considered but, as seen in Section~\ref{sec:est}, exponential models can be fitted very quickly using a tricky analogy with generalized linear models.

As a direct consequence of~\eqref{eq:model2}, the distribution of $\bX$ is necessarily non-stationary. It is therefore highly relevant to ask the following questions: (A) given $\bbeta_\rho$, $\bz_\rho$, are there models with intensity $\rho$? (B) given $\bbeta_\lambda$ and $\bz_\lambda$, are there models with Papangelou conditional intensity $\lambda$. Answer to (A) is easy, as the Poisson point process already answers to this question. It is also quite simple to design  inhomogeneous Cox point process or determinantal point process to achieve this task (see e.g. \cite{choiruddin2018convex,lavancier2015determinantal}). The question for (B) is much more complex (at least if $l\ge 1$, otherwise we are back to the Poisson case). As seen from Table~\ref{tab:models}, the question is essentially related to the existence of non-stationary Gibbs models. \cite{dereudre2012existence} (and the references therein) is one of the most popular existence result in the stationary case. In the non-stationary case, \cite{vasseur:dereudre:20} provides sufficient conditions (which already cover a large class of examples): there exists at least one Gibbs measure with Papangelou conditional intensity $\lambda$ if it satisfies for any $u\in \RR^d$ and $\bx \in \mathcal N$
\begin{equation}\label{eq:lsfr}
\lambda(u,\bx:\bbeta_\lambda)=\lambda\{u,\bx \cap B(u,R) ;\bbeta_\lambda\}
\quad \text{ and }   \quad
\lambda(u,\bx;\bbeta_\lambda) \leq \bar \lambda
\end{equation}
where $R,\bar{\lambda}<\infty$. The first part (finite range property) means that the Papangelou conditional intensity at $u$ depends only on points of $\bx$ close to $u$. The second one, called local stability property, tells that the process is stochastically dominated by a Poisson point process. To set the ideas, the inhomogeneous Strauss model with $l=1$ and $s_1(u,\bx)= \sum_{v\in  \bx}\mathbf{1}(\|v-u\|\le R)$ (number of $R$-closed neighbors of $u$ in $\bx$) satisfies~\eqref{eq:lsfr} for any $\psi \in [0,1]$ and $R<\infty$. We refer to \cite{baddeley2015spatial,ba2023inference} for more complex examples. 

Hence, the problem of inferring $\rho$ or $\lambda$ given by~\eqref{eq:model2} is a well-posed one. The aim of next sections is to estimate $\bbeta_\bullet$ based on a single observation $\bx$ in an observation domain say $D$ of $\bX$, a spatial point process defined on $\RR^d$ with intensity $\rho$ (or conditional intensity $\lambda$).

\section{Statistical inference for low-dimensional parametric $\rho,\lambda$} \label{sec:low}

\subsection{Poisson and pseudo-likelihoods} \label{sec:est}

We consider composite likelihood-based techniques to estimate $\bbeta_\bullet$. In particular we define the ($\log$-)Poisson likelihood and ($\log$-)pseudolikelihood respectively given by
\begin{align} 
\ell_\rho(\bbeta_\rho;\bX) &= \sum_{u\in \bX \cap D} \log \rho(u;\bbeta_\rho) - \int_D \rho(u;\bbeta_\rho) \dd u \label{eq:lrho} \\ 
\ell_\lambda(\bbeta_\lambda;\bX) &= \sum_{u\in \bX \cap D \ominus R} \log \lambda(u,\bX\setminus u; \bbeta_\lambda)- \int_{D\ominus R} \lambda(u,\bX; \bbeta_\lambda) \dd u \label{eq:llambda},
\end{align}
where $D\ominus R$ in \eqref{eq:llambda} stands for the domain eroded by the interaction range $R$. We use the the notation $\ell^{(1)}_\bullet$ and $\ell^{(2)}_\bullet$ to denote the gradient vector and the Hessian matrix  with respect to $\bbeta_\bullet$. In addition, our main result Theorem~\ref{thm} requires the definition of $\bH_\bullet(\bbeta_\bullet;\bX)=-\ell^{(2)}_\bullet(\bX;\bbeta_\bullet)$ which are given by
\begin{align}
\bH_\rho(\bbeta_\rho;\bX) = \bH_\rho(\bbeta_\rho) &= \int_{D} \bz_{\rho}(u) \bz_{\rho}(u)^\top \rho(u;\bbeta_\rho)\dd u \\
\bH_\lambda(\bbeta_\rho;\bX) &= \int_{D\ominus R} \bz_{\lambda}(u,\bX) \bz_{\lambda}(u,\bX)^\top \lambda(u,\bX;\bbeta_\lambda)\dd u
\end{align}
and we also define the matrix $\bH_\lambda(\bbeta_\lambda) = \EE\{\bH_\lambda(\bbeta_\lambda;\bX)\}$.

Let us comment on these methodologies. We focus first on~\eqref{eq:lrho} to estimate $\rho$ given by~\eqref{eq:model2}. When $\bX$ comes from a Poisson point process, as suggested by its name, \eqref{eq:lrho} corresponds to the Poisson likelihood and its maximum was shown to converge by \cite{rathbun1994asymptotic}. 
The method is still an efficient one for general models as the gradient of~\eqref{eq:lrho} is easily seen as an estimating equation by Campbell Theorem~\eqref{eq:Campbell}. \cite{waagepetersen2007estimating} and then \cite{waagepetersen2009two}
proved asymptotic properties for general  point processes under increasing domain asymptotic for large classes of mixing point processes. The efficiency of the method has been improved by \cite{guan2010weighted} and then by~\cite{guan2015quasi} using quasi-likelihood. Infill asymptotic results have been obtained very recently by~\cite{choiruddin2021information}.

Equation~\eqref{eq:llambda} has also an intuitive origin: \cite{jensen1991pseudolikelihood} obtained this criterion as the limit of a certain product of conditional densities. This paper together with \cite{jensen1994asymptotic} were also the first one to establish consistency and asymptotic normality for stationary    exponential family Gibbs models satisfying~\eqref{eq:lsfr} under the increasing domain framework. These results were extended by~\cite{billiot2008maximum,coeurjolly2010asymptotic,dereudre2009campbell} for more general models including infinite range pairwise interaction point processes or non-hereditary Gibbs models. Variant of the pseudo-likelihood include the Takacs-Fiksel method, the logistic regression likelihood and a form of quasi-likelihood respectively studied by~\cite{coeurjolly2012takacs,baddeley2014logistic,coeurjolly2016towards}.
The first results for inhomogeneous models were only proved very recently by~\cite{ba2023inference}, still in the increasing domain framework. The popularity of~\eqref{eq:lrho}-\eqref{eq:llambda} lies without doubt by their simple implementation. This is discussed in the next section.

\subsection{Implementation using the Berman-Turner's approximation}

The main challenge to maximize \eqref{eq:lrho}-\eqref{eq:llambda} is the second term which involves an integral over the observation domain and needs to be approximated numerically. The  Berman-Turner scheme \cite{baddeley2015spatial} is the popular approach conducted by discretizing the integral using both quadrature points and data points by
\begin{align*}
{\int_{D}  \rho(u; \bbeta_\rho) \mathrm{d}u} \approx {\sum_{i=1}^{N+M} w(u_i) \rho(u_i; \bbeta_\rho)}
\quad \text{ and } 
\end{align*}
\begin{align*}
{\int_{D \ominus R}  \lambda(u,\bx; \bbeta_\lambda) \mathrm{d}u} \approx {\sum_{i=1}^{N+M} w(u_i) \lambda(u_i, \bx;\bbeta_{\lambda})},
\end{align*}
where $u_i, i=1,\ldots,N+M$ are quadrature points in $D$ or $D\ominus R$ (depending on the problem) involving $N$ data points and $M$ dummy points and where the $w(u_i)>0$ are quadrature weights such that ${\sum_i w(u_i)}=|D|$ (or $|D\ominus R|$). Using this technique, \eqref{eq:lrho}-\eqref{eq:llambda} are then approximated by
\begin{align}
\label{eq:appx:pois}
\ell_\rho(\bbeta_\rho;\bx) \approx \tilde{\ell_\rho}(\bbeta_\rho;\bx) &= {\sum_{i=1}^{N+M} w_i \{y_i \log \rho_i(\bbeta_\rho) - \rho_i(\bbeta_\rho)\}} \\
\ell_\lambda(\bbeta_\lambda;\bx) \approx \tilde{\ell_\lambda}(\bbeta_\lambda;\bx) &= {\sum_{i=1}^{N+M} w_i \{y_i \log \lambda_i(\bx;\bbeta_{{\lambda}}) - \lambda_i(\bx;\bbeta_{{\lambda}})\}}  \label{eq:appx:pseudo}
\end{align}
where $w_i=w(u_i), y_i=w_i^{-1} \mathbf{1}(u_i \in {\bx \cap D})$ (or $D\ominus R$), $\rho_i(\bbeta_\rho)=\rho (u_i; \bbeta_\rho)$ and $\lambda_i(\bbeta_\lambda)=\lambda (u_i,\bx; \bbeta_\rho)$.
It is now relevant to remark that Equations~\eqref{eq:appx:pois}-\eqref{eq:appx:pseudo} are equivalent to a weighted quasi-likelihood function of independent Poisson variables $y_i$ with weights $w_i$.
Therefore, the implementation can take advantage of any software implementing generalized linear models. These methods are in particular implemented in the \texttt{spatstat} \texttt{R} package \cite{baddeley2015spatial}. The accuracy of the approximation of $\ell_\bullet$ by $\tilde \ell_\bullet$ increases when $N$ is small with respect to $M$. If $N$ is too large or if increasing $M$ leads to numerical instabilities, the induced bias can be non negligible. In these situations, alternatives based on the use of an approximation of~\eqref{eq:lrho}-\eqref{eq:llambda} by a logistic regression likelihood are available (see e.g.~\cite{waagepetersen2007estimating,baddeley2014logistic}).

To sum up Section~\ref{sec:low}: when $p$ is moderate, we have at our disposal a bunch of statistical methodologies to estimate either the intensity or conditional intensity function. These methodologies are well-studied from a mathematical point of view and efficiently implemented.


\section{Inference for sparse (conditional) intensity} \label{sec:high}

\subsection{Setting and additional notation}

The application summarized by Figure~\ref{fig:bei} suggests that more refined methods are necessary. The number of spatial covariates is large (close to 100 if one considers topographic, levels of soil nutrients and levels of combinations of soil nutrients) which leads to numerical problems if one considers methods described in Section~\ref{sec:est}. Known as the curse of dimensionality, these problems can be alleviated if one assumes sparsity in the intensity model and makes use of regularized versions of~\eqref{eq:lrho}-\eqref{eq:llambda}. 

Let us turn to the setting of the present paper. We assume that we observe one realization from $\bX_n$ where $(\bX_n)_{n\ge 1}$ is a sequence of point processes defined in $\RR^d$ and observed in $D_n$. We assume that either the intensity $\rho$ or the Papangelou conditional intensity $\lambda$ is modelled by~\eqref{eq:model2}. The true parameter (to be estimated) is denoted by $\bbeta_{0,\bullet}$ and we assume it can be decomposed as $\bbeta_{0,\rho}=(\bbeta_{01,\rho}^\top,\bbeta_{02,\rho}^\top)^\top$ and  $\bbeta_{0,\lambda}=(\bbeta_{01,\lambda}^\top,\bbeta_{02,\lambda}^\top,\bpsi)^\top$ where $\bbeta_{01,\bullet}=0$ and where all components of $\bbeta_{02,\bullet}$ are non zero. We index any (random) vector, matrix in the same way. Thus $\bz_{01,\bullet}$ corresponds to the set of non-informative covariates while the set $\bz_{02,\bullet}$ represents the set of active features. We assume that $\bbeta_{02,\rho}$ or $(\bbeta_{02,\lambda}^\top,\bpsi^\top)^\top$ has  length $s_n$. Thus $\bbeta_{01,\bullet}$ has length $p_n-s_n$. The sequences $s_n$ and $p_n$ may increase with $n$. Finally, to quantify the amount of increase of data with $n$,  we assume for both problems that $\rho,\lambda$, $\bbeta_{0,\bullet}$ and $D_n$ are such that $\mu_n \to \infty$ as $n\to \infty$
\begin{equation} \label{eq:mun}
\mu_n = \EE\{N(D_n)\} = \int_{D_n} \rho(u;\bbeta_{0,\rho}) \dd u   = \int_{D_n} \EE \left\{\lambda(u,\bX_n;\bbeta_{0,\lambda})\right\} \dd u.
\end{equation}
We assume that for any $n \ge 1$, the model is well-defined. In particular for $\lambda$, this means the sequence of Papangelou conditional intensity functions satisfies~\eqref{eq:lsfr}. Note that $\mu_n$ is a function of $D_n, \bbeta_{01}, \bz_{1,\bullet}(u)$ and $s_n$. We believe this kind of framework is original and quite general. It embraces the well-known frameworks called increasing domain asymptotics and infill asymptotics. For the increasing domain context, $D_n \to \mathbb R^d$ and usually $\bbeta_{02,\dots}$ depends only on $n$ through $s_n$. For the infill asymptotics, $D_n=D$ is assumed to be a bounded domain of $\mathbb R^d$ and usually $(\bz_{2})_1(u)=1$, $(\bbeta_{02,\bullet})_1=\theta_n \to \infty$ as $n\to \infty$. In some sense, the parameter $\mu_n$ plays the role of the sample size in standard inference. To reduce notation in the following, unless it is ambiguous, we do not index $\bX$, $\rho$, $\lambda$, $\bbeta_0$, $\bbeta$, $\bz_\bullet(u), \ell_\bullet$ with $n$. 

When the number of parameters is large, regularization methods allow one to perform both estimation and variable selection simultaneously. When $p_n=p$, \cite{choiruddin2018convex} consider several regularization procedures which consist in adding a convex or non-convex penalty term to~\eqref{eq:lrho}-\eqref{eq:llambda}. A quite similar approach was considered in~\cite{ba2023inference} for Gibbs point processes. To ease the presentation and focus more on the similarities between the two problems of inferring~\eqref{eq:model2}, we only consider the $\ell^1$ regularization which gives rise to the adaptive lasso procedure. The $\ell^1$-regularized versions of~\eqref{eq:lrho}-\eqref{eq:llambda} are given by
\begin{align} \label{eq:Qbullet}
	Q_\bullet( \bbeta;\bX)&= \frac{1}{\mu_n}\ell_\bullet(\bbeta_\bullet;\bX) - {\sum_{j=1}^{p_n} \tau_{n,j}|\beta_{j,\bullet}|} 
\end{align}
where the real numbers $\tau_{n,j}$ are non-negative tuning parameters. The adaptive lasso estimator is then defined by
\begin{align}
\hat \bbeta_{\bullet}= \arg\max_{\bbeta \in \mathbb{R}^{p_n}} Q_\bullet( \bbeta_\bullet;\bX). \label{eq:defest}
\end{align}
When $\tau_{n,j}=0$ for $j=1,\dots,p_n$, the method reduces to the maximum Poisson likelihood or pseudo-likelihood estimator and when $\tau_{n,j}=\tau_n$ to the standard lasso estimator. Note that in the formulation~\eqref{eq:Qbullet}, if we want the model to necessarily have an intercept term and/or if one does not want to regularize the parameter vector corresponding to the interaction term for $\lambda$, we can simply set the corresponding tuning parameters to 0. The choice of $\mu_n$ as a normalization factor in~\eqref{eq:Qbullet} follows the implementation of the adaptive lasso procedure for generalized linear models in the standard software (e.g. \texttt{R} package \texttt{glmnet}~\cite{friedman2010regularization}).

\subsection{Asymptotic results for the adaptive lasso}

Our result relies upon the following conditions:
\begin{enumerate}[($\mathcal C$.1)]
\item For any $n\ge 1$, the intensity or the conditional intensity functions has the log-linear specification given by~\eqref{eq:model2} where $\bbeta_\bullet \in \mathbb R^{p_n}$.  For $\lambda$, we assume that it satisfies~\eqref{eq:lsfr}. \label{C:model}
\item  $(\mu_n)_{n\ge 1}$ is an increasing sequence of real numbers, such that $\mu_n\to \infty$ as ${n} \to \infty$. \label{C:mun}
\item As $n \to \infty$, $\ell^{(1)}_\bullet (\bbeta_{0,\bullet}; \bX ) = O_\P \left( \sqrt{p_n \mu_n}\right)$. \label{C:l1bullet}
\item The matrix $\bH_\rho(\bbeta_{0,\rho})$ or the matrices $\bH_\lambda(\bbeta_{0,\lambda};\bX)$ and $\bH_\lambda(\bbeta_{0,\lambda})$  satisfy
\[
\inf_{n\ge 1} \inf_{\boldsymbol \phi\in \mathbb R^{p_n}, \|\boldsymbol \phi\|=1}  \mu_n^{-1} \boldsymbol{\phi}^\top \bH_\bullet(\bbeta_{0,\bullet}) \boldsymbol{\phi} >0
\quad \text{ and } 
\]
\[
\sup_{\boldsymbol \phi\in \mathbb R^{p_n}, \|\boldsymbol \phi\|=1}   \boldsymbol{\phi}^\top  \left\{ 
\bH_\lambda(\bbeta_{0,\lambda};\bX)-\bH_\lambda(\bbeta_{0,\lambda}
\right\} \boldsymbol{\phi}
=o_\P(\mu_n).
\]
 \label{C:cov}
\item As $n\to \infty$, $p_n^4/\mu_n\to 0$.
\label{C:snpn}
\item For any $c\in \RR$, $\tilde \bbeta_\bullet=\bbeta_{0,\bullet}+c\sqrt{p_n/\mu_n}$ and $j=1,\dots,p_n-s_n$
\begin{align*}
\int_{D_n} \|\bz_\rho(u)\|^3 \rho(u;\tilde \bbeta_\rho)\dd u = O(p_n^{3/2}),
\end{align*}
\[
\int_{D_n \ominus R} \|\bz_\lambda(u,\bX)\|^3 \lambda(u,\bX;\tilde \bbeta_\lambda) \dd u= O_\P(p_n^{3/2}),
\]
\[
\int_{D_n} |(\bz_\rho)_j(u)|\|\bz_\rho(u)\| \rho(u;\tilde \bbeta_\rho)\dd u = O(\sqrt{p_n})
\quad \text{ and }
\]
\[
\int_{D_n\ominus R} |(\bz_\lambda)_j(u,\bX)|\|\bz_\lambda(u,\bX)\| \lambda(u,\bX;\tilde \bbeta_\lambda) \dd u= O_\P(\sqrt{p_n}).
\]
\label{C:zcube}
\item Let  $a_n=\max_{j=p_n-s_n+1,\ldots,{p_n}} \tau_{n,j}$ and $b_n=\min_{j={1},\ldots,p_n-s_n} \tau_{n,j}$. The $\tau_{n,j}$ are allowed to be stochastic and we assume that, as $n \to \infty$
\begin{align*}
a_n \sqrt{\frac{s_n \mu_n }{p_n}}= o_\P(1) 
\qquad\text{ and } \qquad
\frac{1}{b_n} \sqrt{\frac{p_n^2}{\mu_n}} = o_\P(1).
\end{align*}
\label{C:anbn}
\end{enumerate}

Let us discuss these conditions. 
Conditions~\ref{C:model}-\ref{C:mun} have already been explained. It is worth saying that in an attempt to embrace  both problems of estimating $\rho$ or $\lambda$ (and ease the presentation) some conditions may appear useless or too vague. This is the case for \ref{C:l1bullet}, \ref{C:cov} and~\ref{C:zcube}. We invite the reader to refer to \cite{choiruddin2023adaptive,ba2023inference} where conditions on second-order moments (for $\rho$) or second-order Papangelou conditional intensity (for $\lambda$ and in the increasing domain framework) are presented. Essentially, condition~\ref{C:l1bullet}
is obtained by proving that the variance of the score behaves as $p_n \mu_n$ (which corresponds to $n p_n$ for standard GLMs for instance). Condition~\ref{C:cov} shows that the smallest eigenvalue of $\mu_n^{-1}\bH_\bullet(\bbeta_{0,\bullet})$ is positive for any $n$ and could be compared to an assumption of the form $n^{-1}\bX\bX^\top$ tends to a positive definite matrix   for standard linear models with fixed number of covariates (where $\bX$ would stand for the design matrix). Condition~\ref{C:zcube} looks meaningless however, for instance for $\rho$ and $c=0$, these are fulfilled as soon as $\sup_n \sup_i \sup_u |(\bz)_i(u)|<\infty$ (again see~\cite{choiruddin2023adaptive,ba2023inference} for more details). Condition~\ref{C:snpn} is one of the most important and provides a restriction on the number of covariates. Condition~\ref{C:anbn} expresses the compromise to be considered on $s_n,p_n,\mu_n$ and the regularization parameters to ensure consistency and oracle properties. These are quite similar to the corresponding ones considered by~\cite{choiruddin2023adaptive,ba2023inference}. We include here the possibility to have stochastic regularization parameters. Condition~\ref{C:anbn} shows the interest of the adaptive lasso. Indeed for the standard lasso $a_n=b_n$ and both conditions cannot be fulfilled simultaneously even if $s_n$ and $p_n$ do not depend on $n$.

\begin{theorem} \label{thm}
Let $\hat \bbeta_\bullet$ be given by~\eqref{eq:defest}. Assume that  the conditions \ref{C:model}-\ref{C:anbn} hold, then the following properties hold.
\begin{enumerate}[(i)]
\item Consistency: $\hat \bbeta_\bullet$ satisfies
$\hat\bbeta_{\bullet} -\bbeta_{0,\bullet}= O_{\P} ( \sqrt{p_n/\mu_n})$.
\item Sparsity: $\mathrm{P}(\hat \bbeta_{1,\bullet}=0) \to 1$ as $n \to \infty$.
  \end{enumerate}
\end{theorem}

A look at the proof of Theorem~\ref{thm}(i) shows in particular that the consistency remains true under the first part of condition~\ref{C:anbn} and so  remains valid for the standard lasso and actually even if no regularization is considered.  The combination of Theorem~\ref{thm}~(i)-(ii) justifies the interest of a regularization technique. With the same rate of convergence than the one of the unregularized estimator, we ensure that the parameters estimates for non-informative covariates can be set to 0 with probability tending to 1. 

We decided not to include more results to keep the paper readable and short. \cite{choiruddin2023adaptive,ba2023inference} explore for example the asymptotic normality of $\hat \bbeta_{2,\rho}$ and $(\hat \bbeta_{2,\lambda}^\top,\hat \bpsi^\top)^\top$ respectively, and estimates of asymptotic covariance matrices. In the same vein, we only focus on the adaptive lasso penalty. More convex penalties (e.g. adaptive elastic net) as well as non-convex penalties (such as MC+, SCAD penalties) are considered in aforementioned references.

Theorem~\ref{thm} extends the results obtained by~\cite{choiruddin2023adaptive,ba2023inference}. In particular for the estimation problem of $\lambda$, we consider the possibility to include an infill asymptotics and we allow the number of parameters for the interaction parameter, that is the length of $\bpsi$, to increase with $n$. For both problems, we allow the tuning parameters to be stochastic, see Section~\ref{sec:numAL} and this of great interest. Indeed, \cite{zou2009adaptive} and many authors after this work suggest to scale the regularization parameters as $\tau_{n,j}=\tau_n/|(\check \bbeta_{\bullet})_j|^\gamma$ where $\tau_n$ is a non-negative sequence, $\gamma>0$ and where $\check \bbeta_\bullet$ is the unregularized estimator (or actually any estimate producing no zero). 
This popular idea is indeed very natural as it is expected that the $\tau_{n,j}$ for the non-informative covariates will be much larger than the ones for the informative ones and so estimates of parameters corresponding to  non-informative covariates will be more likely set to zero thanks to the $\ell^1$ penalty. We claim that $\tau_n$ can be adjusted to fulfill condition~\ref{C:anbn}. Indeed, since $(\check \bbeta_{\bullet})_j-(\bbeta_{0,\bullet})_j=O_\P(\sqrt{p_n/\mu_n})$, it is easily deduced that $a_n=O_\P(\tau_n)$ and $b_n^{-1}=\max_{j\leq p_n-s_n} \tau_{n,j}^{-1}=O_\P\{\tau_n^{-1} (p_n/\mu_n)^{\gamma/2}\}$. Now, since $s_n\le p_n$ and $p_n^4/\mu_n=o(1)$ by condition~\ref{C:snpn}, 
\begin{equation*}
a_n\sqrt{\frac{s_n\mu_n}{p_n}} = O_\P(\tau_n \sqrt{\mu_n}) 
\quad \text{ and } \quad
\frac1{b_n}\sqrt{\frac{p_n^2}{\mu_n}}=o_\P\left(\frac1{\tau_n} \frac{1}{\mu_n^{-1/4-3\gamma/8}} \right).
\end{equation*}
And so condition~\ref{C:anbn} is in particular satisfied if $\tau_n\sqrt{\mu_n}\to 0$ and $\tau_n \mu_n^{1/4+3\gamma/8} \to \infty$ as $n\to \infty$. For instance, if $\tau_n=\mu_n^{-\alpha}$ with $\alpha>0$, this imposes the non-empty condition $1/2<\alpha<1/4+3\gamma/8$ (if $\gamma>2/3$). A discussion on how to select $\tau_n$ and/or $\gamma$ and other numerical considerations are presented in the next section.

\subsection{Numerical considerations, algorithms and implementation}
\label{sec:numAL}

To implement regularization methods for spatial point processes in particular within \texttt{R}, we combine the existing \texttt{R} package \texttt{spatstat}~\cite{baddeley2015spatial} (devoted to the analysis of spatial point pattern data) with two \texttt{R} packages \texttt{glmnet}~\cite{friedman2010regularization} (for convex penalties) and \texttt{ncvreg}~\cite{breheny2011coordinate} (for non-convex penalties) which use coordinate descent procedures/algorithms~\cite{friedman2010regularization,breheny2011coordinate,choiruddin2018convex,daniel2018penalized}. Those methods rely on the tuning parameters $\tau_{n,j}$. Following \cite{zou2009adaptive}, we suggest to use $\tau_{n,j}=\tau_n/|(\check \bbeta_{\bullet})_j|^\gamma$ where $\tau=\tau_n$ is to be chosen (the parameter $\gamma$ has lesss influence and is often set to 1). Large values of $\tau$ yield estimates with high biases and low variances, whereas small values of $\tau$ produce estimates with low biases and high variances. Therefore, an optimal choice of the tuning parameter $\tau$ is necessary to control the trade-off between the bias and the variance. To select $\tau$, it is reasonable first to identify a decreasing sequence of $\tau$ ranging from a maximum value of $\tau$ for which all penalized coefficients are zero to $\tau=0$ (which corresponds to the unregularized parameter estimates); and second to define a criterion to select $\tau$ by an optimization (minimization) procedure. Following \cite{choiruddin2021information,ba2023inference}, we suggest the use of an information criteria such as the (composite) Bayesian information criterion (BIC) \cite{schwarz1978estimating,gao2010composite} or the (composite) extended regularized information criterion (ERIC) \cite{hui2015tuning} to select $\tau$. Let us define first the (composite) BIC, which we denote by cBIC, 
\begin{align}
\text{cBIC}_\bullet (\tau) = -2 \ell_\bullet(\hat \bbeta_\bullet;\bx) + \log(N) \, d_\bullet(\tau) 
\end{align}
where $N$ is the observed number of points and
\begin{align}
d_\bullet(\tau) = \text{trace}( \hat \bH_\bullet( \hat \bbeta_\bullet) \hat \bSigma_\bullet(\hat \bbeta_\bullet))
\end{align}
with $\bSigma_\bullet(\bbeta_\bullet)  = \bH_\bullet (\bbeta_\bullet)^{-1} \bV_\bullet(\bbeta_\bullet) \bH_\bullet (\bbeta_\bullet)^{-1}$ and $\bV_\bullet(\bbeta_\bullet) = \Var(\ell^{(1)}_\bullet (\bbeta_\bullet;\bX))$. It is worth pointing out that $d(\tau)$ is called the effective number of parameters in the model with tuning parameter $\tau$ and for models with tractable likelihood functions like the inhomogeneous Poisson point process, $d(\tau)$ corresponds to the number of non-zero coefficients in $\hat \bbeta_\bullet$ and the criterion reduces to BIC. For Gibbs models, estimates of $\bH_\lambda$ and $\bSigma_\lambda$ can be efficiently computed using the \texttt{vcov} function of the \texttt{spatstat R} package \cite{coeurjolly2013fast}. Let us now define the (composite) ERIC, which we denote by cERIC and which is designed for the purpose of taking into account the effects of the tuning parameter $\tau$,
\begin{align}
\text{cERIC}_\bullet (\tau) = -2 \ell_\bullet(\hat \bbeta_\bullet;\bx) + \log \left( \frac{N}{|D| \tau} \right)  d_\bullet(\tau).
\end{align}
To sum up, we choose the tuning parameter $\tau \ge 0$ which minimizes either cBIC or cERIC. For more numerical details as well as implementation of the methodology, we refer the reader to \cite{ba2023inference,choiruddin2018convex,choiruddin2021information,choiruddin2023combining}. We end this section by mentioning that a recent version of the \texttt{spatstat R} package allows to include elastic regularization in the \texttt{ppm} function through the option \texttt{improve.type=``enet''}.   



\appendix

\section{Proof of Theorem~\ref{thm}} \label{sec:proof}

\begin{proof} (i)
Let $\bk \in \RR^{p_n}$. We remind the reader that the estimator of $\bbeta_{0,\bullet}$ defined as the maximum of $Q_\bullet$ given by~\eqref{eq:Qbullet}. We aim at proving that for any given $\varepsilon>0$, there exists sufficiently large $K$ such that for 
$n$ sufficiently large 
\begin{align} 
\P \left\{ \sup_{\|\bk\|=K} \Delta_\bullet(\bk;\bX)>0\right\} \le \varepsilon
\label{eq:Delta}
\end{align}
where $\Delta_\bullet(\bk; \bX)= Q_\bullet(\bbeta_{0,\bullet}+\sqrt{p_n/\mu_n}\bk;\bX)-Q_\bullet(\bbeta_{0,\bullet};\bX).$
Equation~\eqref{eq:Delta} will imply that with probability at least $1-\varepsilon$, there exists a local maximum in the ball $\{\bbeta_{0,\bullet}+\sqrt{p_n/\mu_n}\bk: , \|\bk\|\le K\}$. We decompose $\Delta_\bullet(\bk;\bX)=T_{1,\bullet}+T_{2,\bullet}$ with
\begin{align}
 T_{1,\bullet} &= \mu_n^{-1}\left\{
\ell_\bullet(\bbeta_{0,\bullet}+\sqrt{p_n/\mu_n}\bk)-\ell_\bullet(\bbeta_{0,\bullet};\bX)
 \right\} \label{eq:T1}\\
 T_{2,\bullet} & =
 \sum_{j=1}^{p_n} \tau_{n,j} \left( 
 |(\bbeta_{0,\bullet})_j| -|(\bbeta_{0,\bullet})_j+\sqrt{p_n/\mu_n}k_j|
 \right).
 \label{eq:T2}
\end{align}
Since $\rho(u;\cdot)$ and $\lambda(u,\bx;\cdot)$ are infinitely continuously differentiable for any $u\in \RR^d$ and $\bx \in \mathcal N$, $\ell_\bullet(\dot;\bX)$ is in particular twice continuously differentiable. Using a second-order Taylor expansion there exists $t\in (0,1)$ such that
$$  
\mu_n T_{1,\bullet}= \sqrt{\frac{p_n}{\mu_n}} \bk^\top \ell^{(1)}_\bullet (\bbeta_{0,\bullet};\bX) + T_{11,\bullet}+T_{12,\bullet}
$$
where (remind that by definition $\ell^{(2)}_\bullet=-\bH_\bullet$)
\begin{align}
T_{11,\bullet} &= -\frac12 \frac{p_n}{\mu_n} \bk^\top\bH_\bullet(\bbeta_{0,\bullet};\bX) \bk \\
T_{12,\bullet} &= \frac12 \frac{p_n}{\mu_n} \bk^\top
\left\{
\bH_\bullet(\bbeta_{0,\bullet};\bX) -\bH_\bullet(\bbeta_{0,\bullet}+t\sqrt{p_n/\mu_n};\bX) 
\right\}
\bk     .
\end{align}
Under condition~\ref{C:cov}, $T_{11,\rho}\le -(\alpha_\rho/2)p_n\|\bk\|^2$ for some $\alpha_\rho>0$. For the conditional intensity, again using condition~\ref{C:cov}
\begin{align*}
T_{11,\lambda} &= -\frac12 \frac{p_n}{\mu_n} \bk^\top \bH_\lambda(\bbeta_{0,\lambda})\bk + \omega_{11,\lambda} \le  -\frac{\alpha_\lambda}2 p_n\|\bk\|^2
 + \omega_{11,\lambda}
\end{align*}
where $\omega_{11,\lambda}=o_\P(p_n)$. Now, for some $\tilde \bbeta_\bullet$ on the line segment between $\bbeta_{0,\bullet}$ and $\bbeta_{0,\bullet}+t\sqrt{p_n/\mu_n}$, we have under condition~\ref{C:zcube}
\begin{align*}
T_{12,\rho} & =\frac12 \frac{p_n}{\mu_n} \bk^\top \left\{
\int_{D_n} \bz_{\rho(u)}\bz_\rho(u)^\top t\sqrt{\frac{p_n}{\mu_n}} \bk^\top \bz_\rho(u) \rho(u;\tilde \bbeta_{\rho} \dd u
\right\} \bk\\
&= O \left(\frac{p_n}{\mu_n} \sqrt{\frac{p_n}{\mu_n}}\right) \int_{D_n} \|\bz_\rho(u)\|^3 \rho(u;\tilde \bbeta_\rho)\dd u = O\left( p_n \sqrt{\frac{p_n^4}{\mu_n}}\right)\\
T_{12,\lambda} & =\frac12 \frac{p_n}{\mu_n} \bk^\top \left\{
\int_{D_n\ominus R} \bz_{\lambda(u,\bX)}\bz_\lambda(u,\bX)^\top t\sqrt{\frac{p_n}{\mu_n}} \bk^\top \bz_\lambda(u,\bX) \lambda(u,\bX;\tilde \bbeta_{\lambda} \dd u
\right\} \bk\\
&= O \left(\frac{p_n}{\mu_n} \sqrt{\frac{p_n}{\mu_n}}\right) \int_{D_n\ominus R} \|\bz_\lambda(u,\bX)\|^3 \lambda(u,\bX;\tilde \bbeta_\lambda)\dd u = O_\P\left( p_n \sqrt{\frac{p_n^4}{\mu_n}}\right).
\end{align*}
Hence, under condition \ref{C:snpn}, $T_{12,\rho}=o(p_n)$ and $T_{12,\lambda}=o_\P(p_n)$, which yields
\begin{align*}
 T_{1,\bullet} \le \frac{1}{\mu_n} \sqrt{\frac{p_n}{\mu_n}} \bk^\top \ell^{(1)}_\bullet(\bbeta_{0,\bullet}) - \frac{\alpha_\bullet}2 \frac{p_n}{\mu_n} \|\bk\|^2 + \omega_{1,\bullet}
 \end{align*}
 where $\omega_{1,\rho}=o(p_n/\mu_n)$ and $\omega_{1,\lambda}=o_\P(p_n/\mu_n)$. Regarding the term $T_{2,\bullet}$, we have
 \begin{align*}
 T_{2,\bullet} & = \sum_{j=1}^{p_n-s_n} \tau_{n,j} \left( 
 |(\bbeta_{0,\bullet})_j| -|(\bbeta_{0,\bullet})_j+\sqrt{p_n/\mu_n}k_j|
 \right)
+ \\
&\quad \quad \quad  \quad \quad \sum_{j=p_n-s_n+1}^{p_n} \tau_{n,j} \left( 
 |(\bbeta_{0,\bullet})_j| -|(\bbeta_{0,\bullet})_j+\sqrt{p_n/\mu_n}k_j|
 \right)\\
 &\le\sum_{j=p_n-s_n+1}^{p_n} \tau_{n,j} \left( 
 |(\bbeta_{0,\bullet})_j| -|(\bbeta_{0,\bullet})_j+\sqrt{p_n/\mu_n}k_j|
 \right)\\
 &\le a_n\sqrt{\frac{p_n}{\mu_n}} \sum_{j=p_n-s_n+1}^{p_n} |k_j|\le a_n \sqrt{\frac{s_np_n}{\mu_n}} \|\bk\| =: \omega_{2,\bullet}
 \end{align*}
 where under condition~\ref{C:anbn}, $\omega_{2,\bullet}=O(a_n\sqrt{{s_np_n}/{\mu_n}} {p_n}/{\mu_n})= o_\P(p_n/\mu_n)$. We deduce that as $n\to \infty$
 \begin{align*}
\Delta_\bullet(\bk;\bX) \le \frac1\mu_n \sqrt{\frac{p_n}{\mu_n}} \|\ell^{(1)}_\bullet(\bbeta_{0,\bullet};\bX)\| \, \|\bk\| -\frac{\alpha_\bullet}2 \frac{p_n}{\mu_n} \|\bk\|^2 + \omega_\bullet   \end{align*}
where $\omega_\bullet = \omega_{1,\bullet}+\omega_{2,\bullet}=o_\P(p_n)$ whereby we continue with
\begin{align*}
\P \left\{ \sup_{\|\bk\|=K} \Delta_\bullet(\bk;\bX)>0\right\} & 
\le \P \left( L_n \ge \frac{\alpha_\bullet}2 K \sqrt{p_n\mu_n}
\right)    
\end{align*}
where $L_n=\|\ell^{(1)}_\bullet(\bbeta_{0,\bullet};\bX)\| + \mu_n\sqrt{\frac{\mu_n}{p_n}} \omega_\bullet$. This leads to the result since under condition~\ref{C:l1bullet}
\begin{align*}
L_n=O_P\left(\sqrt{p_n \mu_n}\right) + O\left(\mu_n\sqrt{\frac{\mu_n}{p_n}} \omega_\bullet\right)    = O_P\left(\sqrt{p_n \mu_n}\right).
\end{align*}

(ii) Following (i), we intend to prove that for any $\tilde \bbeta_\bullet\in \RR^{s_n}$ (where $\tilde \bbeta_\rho=\bbeta_{2,\rho}$ or $\tilde \bbeta_\lambda=(\bbeta_{2,\lambda}^\top,\bpsi^\top)^\top$) satisfying $\|\tilde \bbeta_\bullet-\tilde\bbeta_{0,\bullet}\|=O_\P (\sqrt{p_n/\mu_n})$ and any $K_1>0$
\begin{align*}
   Q_\bullet \left\{ (\mathbf{0}^\top, \tilde \bbeta_\bullet^\top)^\top; \bX\right\} = \max_{\|\bbeta_{1,\bullet}\|\le K_1\sqrt{p_n/\mu_n}} Q_\bullet \left\{ (\bbeta_{1,\bullet}^\top, \tilde \bbeta_\bullet^\top)^\top; \bX\right\} .
\end{align*}
To this end, it is sufficient to show that with probability tending to 1 as $n\to \infty$, for any $\tilde\bbeta_{\bullet}$ such that $\|\tilde \bbeta_\bullet-\tilde \bbeta_{0,\bullet}\|=O_\P(\sqrt{p_n/\mu_n})$, we have for any $j=1,\dots,p_n-s_n$
\begin{align}
\frac{\partial Q_\bullet(\bbeta_\bullet;\bX)}{\partial(\bbeta_\bullet)_j} <0 
 \; \text{ for } \;  &
0<(\bbeta_\bullet)_j<\varepsilon_n, \;    
\frac{\partial Q_\bullet(\bbeta_\bullet;\bX)}{\partial(\bbeta_\bullet)_j} >0 
\; \text{ for } \;
-\varepsilon_n<(\bbeta_\bullet)_j<0.\label{eq:derivQ}
\end{align}
We focus only on the first part of~\eqref{eq:derivQ} as the other one follows along similar lines. We have
\begin{align*}
\frac{\partial \ell_\bullet(\bbeta_\bullet;\bX)}{\partial(\bbeta_\bullet)_j}  =
\frac{\partial \ell_\bullet(\bbeta_{0,\bullet};\bX)}{\partial(\bbeta_\bullet)_j}  +
R_\bullet
\end{align*}
where
\begin{align*}
R_\rho &= -\int_{D_n} (\bz_\rho)_j(u) 
\left\{\rho(u;\bbeta_\rho)-\rho(u;\bbeta_{0,\rho}) \right\} \dd u\\
R_\lambda &= -\int_{D_n} (\bz_\lambda)_j(u,\bX) 
\left\{\lambda(u,\bX;\bbeta_\lambda)-\lambda(u,\bX;\bbeta_{0,\lambda}) \right\} \dd u.
\end{align*}
By Taylor expansion and Cauchy-Schwarz inequality, there exists $\check \bbeta_\bullet$ on the line segment between $\bbeta_{0,\bullet}$ and $\bbeta_\bullet$
\begin{align*}
|R_\rho| &= O\left(\|\bbeta_\rho-\bbeta_{0,\rho}\|\right) \int_{D_n} |(\bz_\rho)_j(u)| \, \|\bz_\rho(u)\| \rho(u;\check\bbeta_\rho) \dd u    \\
|R_\lambda| &= O\left(\|\bbeta_\lambda-\bbeta_{0,\lambda}\|\right) \int_{D_n\ominus R} |(\bz_\lambda)_j(u,\bX)| \, \|\bz_\lambda(u,\bX)\| \lambda(u,\bX;\check\bbeta_\lambda) \dd u.    
\end{align*}
By condition \ref{C:zcube}, we deduce that $R_\bullet = O_\P(\sqrt{p_n/\mu_n} \sqrt{p_n}\mu_n)=O_\P(p_n\sqrt{\mu_n})$, which combined with condition~\ref{C:l1bullet} yields
\begin{align}
    \frac{\partial \ell_\bullet(\bbeta_\bullet;\bX)}{\partial(\bbeta_\bullet)_j}  =O_\P(p_n\sqrt{\mu_n}).\label{eq:derivell}
\end{align}
Now, let $0<(\bbeta_\bullet)_j<\varepsilon_n$, for $n$ sufficiently large
\begin{align*}
\P \left\{ 
\frac{\partial Q_\bullet(\bbeta_\bullet;\bX)}{\partial(\bbeta_\bullet)_j} <0
\right\} &=\P \left\{
\frac{\partial \ell_\bullet(\bbeta_\bullet;\bX)}{\partial(\bbeta_\bullet)_j} - \mu_n \tau_{n,j} \mathrm{sign}(\bbeta_\bullet)_j <0
\right\} \\
&=\P \left\{
\frac{\partial \ell_\bullet(\bbeta_\bullet;\bX)}{\partial(\bbeta_\bullet)_j} < \mu_n \tau_{n,j} 
\right\} \\
&\ge \P \left\{
\frac1{b_n} \, \frac{\partial \ell_\bullet(\bbeta_\bullet;\bX)}{\partial(\bbeta_\bullet)_j} <\mu_n
\right\}
\end{align*}
which tends to 1 as $n\to \infty$ since by condition~\ref{C:anbn} and~\eqref{eq:derivell}
\begin{align*}
    \frac1{b_n} \, \frac{\partial \ell_\bullet(\bbeta_\bullet;\bX)}{\partial(\bbeta_\bullet)_j} = o_\P\left( \sqrt{\frac{\mu_n}{p_n}}\right) \, O_\P\left(p_n \sqrt{\mu_n} \right) = o_\P(\mu_n).
\end{align*}
\end{proof}

\bibliographystyle{imsart-nameyear}
\bibliography{main}

\end{document}